\documentclass[a4paper,10pt]{article}
\usepackage{amsfonts}
\usepackage{amssymb}
\usepackage[utf8]{inputenc}

\usepackage[english ]{babel}

\usepackage{bbding}

\usepackage{MnSymbol}

\usepackage[all]{xy}

\usepackage{tikz-cd} 
\usepackage[all]{xy}
\usepackage{amsfonts}
\usepackage{amssymb}
\usepackage{graphicx}
\usepackage{mathtools}
\usepackage{skak} 
\usepackage{foekfont}

\usepackage{calligra}
\usepackage[T1]{fontenc}

\usepackage{amsmath,mathdots}

\usepackage{stmaryrd}

\usepackage[utf8]{inputenc}

\usepackage{mathtools}

\usepackage{tikz}

\definecolor{amber(sae/ece)}{rgb}{1.0, 0.49, 0.0}
\newfont{\rsfsten}{rsfs10 scaled 1200}
\usepackage{MnSymbol}
\usepackage{tikz}
\usepackage{amscd}
\usepackage{MnSymbol}
\usepackage{wasysym}

\usepackage{mathrsfs}

\usepackage{pdfcolmk}

\usepackage{graphicx}

\newcommand*{\rom}[1]{\expandafter\@slowromancap\romannumeral #1@}

\usepackage{wrapfig}

\newcommand{\tightunderset}[2]{%
  \mathop{#2}\limits_{\vbox to .3ex{\kern-0.95ex\hbox{$#1$}\vss}}}
\newcommand{\tightoverset}[2]
{%
  \mathop{#2}\limits_{\vbox to .3ex{\kern-0.95ex\hbox{$#1$}\vss}}}

\newcommand{\oset}[2]{%
  {\mathop{#2}\limits^{\vbox to -.5\ex@{\kern-\tw@\ex@
   \hbox{\scriptsize #1}\vss}}}}
\usepackage[utf8]{inputenc}

\usepackage{stmaryrd}

\usepackage {ifsym}

\usepackage {url}

\usepackage{amsmath,amssymb,amsthm}

\title{ Scalar Curvature of  Manifolds with Boundaries: Natural Questions and Artificial Constructions}

\author{Misha Gromov}

\begin{document}
\maketitle

 During  the workshop {\sl Emerging Topics: Scalar Curvature and Convergence} 
in October 15-19, 2018 at IAS in Princeton, Chao Li,  Pengzi Miao, Andr\'e Neves, Christina Sormani and myself  discussed
a  possible formulation of and  approach to  the solution of  the following    problems. \vspace{1mm}

{\large {\sf Problem} {\sf \textbf  A.}} Let $Y=(Y,h)$ be a closed $(n-1)$-dimensional  Riemannin manifold, which, besides a Riemannian metric $h$,   carries  a continuous  function $M=M(y)$ on it.

Find condition(s) on $(Y,h, M) $   that would allow/wouldn't allow \vspace{1mm}

{\sl  a filling of     $Y$ by a  compact Riemannin $(n+1)$-dimensional  manifold $X =(X, g)$, where "filling" means that
$$\partial X=Y, $$ 
where

$\bullet$ the restriction of the Riemannin metric $g$ to $Y$ is equal to $h$;
   
   $\bullet$ the mean curvature of $Y$ in $X$ is equal to $M$,}

 (by our sign convention  {\it convex} boundaries have    $M\geq 0$) {\sl   
  and where the essential property required of $X$ is
 
     $\bullet$ non-negativity of the scalar curvature,  
   $$Sc(X)=Sc(g)\geq 0,$$
   or, more generally, a lower bound $Sc(X)\geq \sigma$ for a given $ \sigma\in (-\infty,  \infty)$. }
   \vspace{1mm}

  {\large {\sf Problem} \sf \textbf B.} Granted that  a filling $X$ with $Sc(X)\geq 0$ (or with $Sc(X)\geq \sigma$) exists, what are the constraints  on the geometry of such an  $X$ imposed by
   $(Y,h, M)$? \vspace{1mm}

  For example: \vspace{1mm}

 {\sf \textbf A$_1$}. Does  {\it sufficient}, depending on $(Y,g)$, {\it mean convexity}, i.e. "large  positivity" of the mean curvature of $Y$,  rule out such fillings? 
   
  {\sf \textbf B$_1$}. Is there a {\it lower bound on the volume}  of  filling  manifolds $X $,   
   in terms  of $(Y,h,M)$ and the lower bound on the scalar curvature of $X$ (with a particular attention to the  the case $Sc(X)\geq 0)$?\footnote{This question was raised by Pengzi Miao, as I recall.}
   
   \vspace {1mm}
   
   {\sf \textbf C. }  Is there a  relation(s) between {  \sf \textbf A$_1$} and {\it the  lower bounds on the dihedral angles} of Riemannin manifolds $X$ {\it with corners}, where   \vspace {1mm}
   
   \hspace {27mm} $Sc(X)\geq 0$ and 
$mean.curv(\partial X)\geq 0$?\footnote{This question  by    Christina Sormani started our conversation.}
    
    \vspace {1mm}
    
    What follows is (a rectified version  of the first draft)   of what came out  of  our conversation.

    \tableofcontents
   
   \section{Bound on the Size of $\partial X$ by the Scalar Curvature and the Mean Curvature}

    The following inequalities  
   $[mean]_{Sc< 0}$  and  $[mean]_{Sc\geq 0}$, 
   deliver a positive answer to  {\sf \textbf A$_1$} for {\it spin } manifolds $X$.    
     \vspace{1mm}


  {\sf Let $Rad_{S^{n-1}}(Y)$, denote the the {\it (hyper)spherical radius}  of a closed  orientable Riemannin $(n-1)$-dimensional manifold   $Y$, that is the maximal radius  $R$  of the  sphere, such that $Y$ admits a   distance decreasing map $Y\to S^{n-1}(R)$ of non-zero degree. 
 
Observe that\vspace {1mm}

$\bullet $ \hspace {1mm}  $Rad_{S^{n-1+m}}(Y\times S^m(R))=\min(R, Rad_{S^{n-1}}(Y))$;

$\bullet $ \hspace {1mm} if we   require  our {\it distance decreasing}  map  to be  {\it smooth}, we come up with the same  (which could be, a  priori, smaller) value of  $Rad_{S^{n-1}}(Y)$.}

\vspace {1mm}\vspace {1mm}

{\it \large Theorem.} {\sl Let  $X$ be a compact   orientable   spin manifold of dimension $n$ with boundary $Y=\partial X$.

 If $Sc(X)\geq -\sigma $, $\sigma\geq 0$,   then the infimum of the   mean curvature 
of $Y\subset X$ is bounded by
$$\inf_{y\in Y}mean.curv (Y,y)  \leq \max\left (\frac {n+m-1}{Rad_{S^{n-1}}(Y)},  \sqrt {\frac {\sigma}{m(m-1)}}\right). \leqno {[mean]_{Sc< 0}} $$}
for all $m\geq 2$.

Furthermore if $Sc(X)\geq 0$, then

$$\inf_{y\in Y}mean.curv (Y,y) \leq  \frac {n-1}{ Rad_{S^{n-1}}(Y)}. \leqno {[mean]_{Sc\geq 0}}$$


{\it About the Proof.}  Inequality $[mean]_{Sc\geq 0}$ is proven in the next section by mapping  the  double $X\cup_YX$ to the sphere $S^n=B^n\cup_{S^{n-1}}B^n$ endowed  with a radial   metric, which has {\it positive  curvature operator} and  to which 
{\sl  Goette-Semmelmann's  extremality theorem}  [GS 2000] applies.

Then  inequality $[mean]_{Sc< 0}$ follows by applying   $[mean]_{Sc\geq 0}$  to $X\times S^m(R)$
 for $R=\sqrt{\frac{m(m-1)}{\sigma}}$. 
   \vspace{1mm}

{\it Remarks, Examples, Questions.} (a) The inequalities $[mean]_{Sc\geq 0}$ and  $[mean]_{Sc< 0}$
are already significant for closed {\it imbedded} hypersurfaces $Y$ in the Euclidean space $\mathbb R^{n}$ (where they bound domains $X$ with Sc(X)=0) and in the hyperbolic space $\mathbb H^n_{-1}$ (where $Sc(X)=-n(n-1)$).

There is no  apparent {\it direct}  proof of this  with {\it no    Dirac operators behind the scenes}, which are essential   for the proof of Goette-Semmelmann's  theorem.   
 
Notice in this regard  that a proof of such an inequality  can't localise near $Y$. In fact,

{\sf {\it immersed  hypersurfaces  with self-intersection}  in $\mathbb R^{n}$   e.g. finite coverings of a  $2$-torus in $\mathbb R^3$,   with $mean.curv >0$,
may have {\it arbitrarily large (hyper)spherical radii} regardless of their mean  curvatures.}


 \vspace{1mm}

(b) The inequality  $[mean]_{Sc\geq 0}$ is sharp, where balls  $X=B^n(R)\subset \mathbb R^n$ with boundaries $Y=S^{n-1}(R)=\partial B^n(R)$ are   {\it extremal}  for this inequality, since 
$[mean]_{Sc\geq 0}$ turns to equality for balls:
$$mean.curv (S^{n-1}(R)) =  \frac {n-1}{R= Rad_{S^{n-1}}(S^{n-1}(R))}.$$

In fact, {\it \large conjecturally,} balls are {\it  rigid}:

{\sf if $Sc(X)\geq 0$, then the equality  $mean.curv (Y) =  \frac {n-1}{ Rad_{S^{n-1}}(Y)}$ must imply that 
$X$ is isometric to a ball  $B^n(R)\subset \mathbb R^n$.}

But our proof  doesn't yield  rigidity due to an  approximation  error in our  argument, when    we smooth the metric in 
the double $X\cup_YX$.
 \vspace{1mm}

Yet, there is a case, where the proof of the rigidity causes no problem, namely, when the manifold $X$ near its boundary $\partial X$ is isometric to a neighbourhood of  a  sphere  in the Euclidean space $\mathbb R^n$.
In fact, our theorem says  in this case that  \vspace {1mm}
 
 {\it if a connected Riemannin spin  manifold $X^+$  with $Sc(X^+)\geq 0$ is isometric at infinity, i.e. outside a compact subset in $X^+$,  to the complement  of an $R$-ball in the Euclidean space $\mathbb R^n$, then   $X^+$ is isometric to $\mathbb R^n$}.  \vspace {1mm}

(Probably,   Goette-Semmelmann's  Dirac theoretic argument behind our proof converges, when $R\to \infty$,  to that by Witten in his proof of {\it the positive mass theorem}.\footnote {The positive mass theorem, which  claims the    rigidity of $\mathbb R^n$ under a class of perturbations which may have {\it non-compact supports},   was proven in 1979 by  Schoen
and Yau   for manifolds of dimensions $\leq7$ and  in 1981  by Witten for  spin manifolds of all dimensions, while the reduction  to the, a priori more restrictive,  Euclidean rigidity  was found by 
  J. Lohkamp in 1999.
 
  The Euclidean rigidity also is an 
  obvious corollary of 
   
   {\it non-existence of  non-flat metrics with $Sc\geq 0$  on the torus and on overtorical manifolds,} 
   
   \hspace {-4mm} that 
 is now established (see [SL 2017] and [L  2018]  for, all possibly non-spin,  manifolds of all dimensions $n$.})

  \vspace{1mm}
 
 (c) Unlike $[mean]_{Sc\geq 0}$, the inequality $[mean]_{Sc< 0}$ is non-sharp.  \vspace{1mm}
 
 {\large \it Conjecturally},   $\inf_{y\in Y}mean.curv (Y,y)$ must be bounded by the mean curvature of the ball of radius 
 $R=Rad_{S^{n-1}}(Y)$ in the hyperbolic $n$-space with constant  sectional curvature 
 $-\sigma/n(n-1)$.

 \vspace{1mm}

(d)  The spin condition for $X$  can be  relaxed   to  {\it spin of the universal covering $\tilde X$} (see remark (d) following formulation of the Goette-Semmelmann theorem in the next section and also section 10 in 
[G 2018]), but it remains {\large \it  problematic} if this condition is necessary at all. 

Yet,  we indicate in section 6 an  approach, although a very artificial one,   to  such    (non-sharp) inequalities  (with extra assumptions on the geometry of $Y$)   by the techniques of  
{\it minimal hypersurfaces}  where spin is irrelevant.

 \vspace{1mm}


  
 
 (e) One can't, in general,  substitute $\inf_ymean.curv(Y,y)$ in the above inequalities by any  integral  characteristic of the function $mean.curv(Y,y)$ but this may be possible  with a  more subtle (spectral?)  geometric invariant  of $Y$ than $Rad_{S^{n-1}}.$  (Compare with {\huge $\star $} {\huge $\star $} in  section 3.)\vspace{1mm}\vspace{1mm}
 
 \vspace{1mm}

 (f)  What are   simple (non-simple?) examples of  {\it extremal/rigid} Riemannian  $n$-manifolds $X$  with boundaries, where you can't  simultaneously increase 
 $Sc(X)$, $mean.curv(\partial X)$ and   $Rad_{S^{n-1}}(\partial X)$?
 
 For instance what is the {\it sharp} bound on the mean curvature of $Y=\partial X$  by 
 $Rad_{S^{n-1}}(Y)$ for manifolds $X$ with $Sc(X)\geq n(n-1)$?
 
(The subtlety of this extremality/rigidity problem for balls in $S^n$ is  revealed by 
 the results by  Brendle, Marques  and Neves [BMN 2010]  and by Brendle and  Marques [BM 2011]:  \vspace {1mm}
 
{\sf  the   existence/nonexistence  of a  {\it deformation of the  spherical  metric} in the interior of such a  ball with {\it increase of the scalar curvature}  depends of the {\it radius} of the ball.})

 \section{Extremal Metrics on Spheres, Doubles and the Proof  of  Inequality $[mean]_{Sc>0}$  }

  Let $ X$ be a  compact Riemannin $n$-manifold with boundary $Y= \partial X$ and 
 let $X_\varepsilon X$ be the  {\it $[^\varepsilon\leftrightarrows_\varepsilon]^{\hspace {-1mm}\Rightcircle}$-rounding} of 
 the double $X\cup_YX$ of  of $ X$  (compare  p. 227 in   
  [GL 1980] and  section 11.4 in 
 [G 2018]) 
that is the boundary of the  $\varepsilon$-neighbourhood of  
$ X_{-\varepsilon}= X_{-\varepsilon}  \times \{0\}\subset X\times \mathbb R$, 
$$ X_\varepsilon  X=\partial U_\varepsilon (X_{-\varepsilon})\subset  X\times \mathbb R,$$ 
where  $ X_{-\varepsilon} \subset X$ is the complement of the $\varepsilon$-neighbourhood of the boundary
 $\partial X\subset  X$.
\vspace {1mm}
 
The hypersurface   $ X_\varepsilon  X\subset  X\times \mathbb R$
 consists of two $2\varepsilon$-equidistantly parallel copies of    $ X_{-\varepsilon} $ and 
a semicircular part  
 ${\hspace {0mm} {\small \Rightcircle_\varepsilon}}= {\hspace {0mm} {\small \Rightcircle_\varepsilon}}( X)=  X_{-\varepsilon}\times S^1_+(\varepsilon)$, that is  a half of the boundary of  the  $\varepsilon$-neighbourhood of $ X_{-\varepsilon} \times \mathbb R$,
as depicted in  figure 8 on p.  227 in  [GL 1980]).

 The principal curvatures $\lambda_1,...,\lambda_n$ of ${\small \Rightcircle_\varepsilon}$ are evaluated, for small $\varepsilon \to 0$,  in terms of the principal curvatures  $\mu_1,..., \mu_{n-1}$ of the boundary $\partial   X$  as follows.
 $$\lambda_i=(\mu_i +O(\varepsilon)) \cdot\cos \theta+o(\varepsilon) \mbox { for  $i=1,..., n-1$ and $\lambda_n=\varepsilon^{-1}+O(1)$,}$$
where $-\frac {\pi}{2}\leq \theta \leq  \frac {\pi}{2}$ denotes  the angular parameter of the (right) semicircle  
$S_+^1(\varepsilon)$
and  the scalar curvature of ${\hspace {0mm} {\small \Rightcircle_\varepsilon}}$, which is expressed by the Gauss theorema egregium, satisfies
$$Sc({\hspace {-1mm} {\small \Rightcircle_\varepsilon}})(\underline x,\theta) = Sc( X)(\underline x)+ (2\varepsilon ^{-1} m(\underline x)+O(1))\cdot \cos \theta+o(1),$$
where $m(v)= mean.curv(\partial  X)(\underline x)$, 
$ x\in \partial  X$, and $-\frac {\pi}{2}\leq \theta \leq  \frac {\pi}{2}$.

{\it Smoothing Remark.} The submanifold  $X_\varepsilon X \subset X\times \mathbb R$ is only $C^1$-smooth and the induced metric in it is only continuous. But it   can be  smoothed as  in section 11.1 of 
 [G 2018] with an arbitrary small perturbation of the metric  of this submanifold and   of its scalar curvature.

\vspace {1mm}

Let  $f$ be a smooth map from the boundary $Y=\partial X$ to the $R$-sphere  $ S^{n-1}(R)\subset \mathbb R^n$  and let $F:X\to \underline X= B^n(R)\subset \mathbb R^n$ be a smooth extension of $f$ to a map from $X$ to the ball bounded by  $S^{n-1}(R)$.

 Observe that, for small $\varepsilon$, the map  $F$ naturally "suspends" to a map between the  "round doubles"  of these manifolds, denoted 
 $$F_{\hspace{-1mm}\small\Rightcircle_\varepsilon}: X_\varepsilon X \to \underline X _\varepsilon\underline X.$$

Clearly, the  differential of $F_{\hspace{-1mm}\small\Rightcircle_\varepsilon}$  away 
from $\small\Rightcircle_{\varepsilon}(X)$ is the same as $dF$,  while  the norm of $F_{\hspace{-1mm}\small\Rightcircle_\varepsilon}$ at $ ( x, \theta) \in
\small\Rightcircle_{\varepsilon}(X)$ satisfies
$$||dF_{\hspace{-1mm}\small\Rightcircle_\varepsilon}( x, \theta)||\leq ||df(x)|| (1+O(\varepsilon)).$$

Now, these  estimates for $X_\varepsilon X$  along with similar ones for $\underline X_\varepsilon \underline X$, where  $\underline X=B^n(R)$, and with the above smoothing  remark allows 
  reduction of the  desired inequality  $[mean]_{Sc\geq 0}$, written as 
$$Rad_{S^{n-1}}(\partial (X)) \leq  \frac {n-1}{\inf mean. curv (\partial X)}.$$
to the following,

\vspace{1mm}

{\large \it Goette-Semmelmann  Extremality Theorem} 
[GS 2000].  {\sf Let $\underline Z= (S^{n}, \underline g)$ be the  $n$-sphere with a Riemannin  metric  $ \underline g$ which has a {\it non-negative curvature operator}, e.g. $\underline Z \subset \mathbb R^{n+1}$
is a smooth convex  hypersurface obtained by  smoothing  $\underline X_\varepsilon \underline X_\varepsilon\subset  \mathbb R^{n+1}$ for $\underline X_\varepsilon =B^n \subset \mathbb R^n$,  with the induced Riemannin metric. \footnote{Manifolds with positive sectional curvatures   isometrically  embeddable to $\mathbb R^{n+2}$, have positive  curvature operators  by     Weinstein's  theorem  
 [W 1970], but for  $X\subset \mathbb R^{n+1}$ this, probably, was   known prior to Weinstein's paper.}

Let $X$ be a closed connected orientable $n$-manifold and $f:X\to\underline X$ be a smooth map. 

{\it If 

$\bullet $ \hspace {1mm} $X$  is spin and the map $f$ has non-zero degree;

$\bullet $ \hspace {1mm}   $||df(x)||\leq 1$ at all points $x\in X$, where $Sc(\underline X,f(x))>0$;

$\bullet $ \hspace {1mm} 
$Sc(X,x)\geq Sc (\underline X,f(x)).$
for all $x\in X$, }
\vspace{1mm}

then
$$Sc(X,x)= Sc (\underline X,f(x))$$
for  all $x\in X$.}\vspace{1mm}

{\it Thus,  the  inequality  $[mean]_{Sc\geq 0}$, hence, $[mean]_{Sc< 0}$,  is established.}

\vspace{2mm}

 {\it Remarks.} (a)  There are  further (less straightforward) applications of  the  Goette-Semmelmann   theorem   to manifolds with boundaries but  we shall  discuss these somewhere else.

 (b) A non-sharp version of the inequality  $[mean]_{Sc\geq 0}$
  can be derived from a  {\it rough  hypersphericity inequality} 
  [GL 1980], the proof of which, unlike  the  Goette-Semmelmann's   theorem. doesn't depend on a  non-trivial curvature computation.

All we actually  need  is the lower bound 
$Sc(X_\varepsilon X)\geq\frac  {1}{3\varepsilon}$  in the $\frac {\varepsilon}{2}$-neighbourhood of  
$\partial X\subset X_\varepsilon X$  for a small $\varepsilon>0$, where we observe that
\vspace {1mm}

{\sl  smooth maps   $f: \partial X\to S^{n-1}$, $n=dim(X)$,  extend to smooth maps $$F_\varepsilon:X_\varepsilon X \to    S^{n}\supset S^{n-1},$$
such that the two components of the complement to the normal  $\frac {\varepsilon}{2}$-neighbourhood
   $\partial X\times \left[-\frac {\varepsilon}{2}, \frac {\varepsilon}{2}\right]\subset  X_\varepsilon X$ of 
$\partial X\subset X_\varepsilon X$ collapse to the two  poles of $S^{n}$
and such that the second exterior power of the differential of $F_\varepsilon$  is bounded by the norm of the differential of $f$ on $\partial X$ as follows. }
$$ ||\wedge^2(dF_\varepsilon(x, \delta)||\leq \frac {10}{\varepsilon}||df_{|\partial X}(x)||  $$
for  all points $(x, \delta)\in \partial X\times \left[-\frac {\varepsilon}{2}, \frac {\varepsilon}{2}\right]$, $x\in 
\partial X$, $\delta \in \left[-\frac {\varepsilon}{2}, \frac {\varepsilon}{2}\right].$

Since $deg(F_\varepsilon)=deg (f)$, the rough  bound
$$Sc( X_\varepsilon X, (x,\delta) \leq const_n  ||\wedge^2(dF_\varepsilon(x, \delta)||\cdot Sc(S^n)$$
 from 
 [GL 1980] yields the corresponding bound 
$$mean.curv (\partial X,x)\leq const'_n ||df_{\partial X}||\cdot mean.curv (S^{n-1}).$$

Notice, that one could apply here    {\it  Llarull's inequality} 
 [Ll 1998]   to $F_\varepsilon$  with {\it sharp}  $const_n=1 $, but this would not make the resulting bound for $f$ sharp anyway.

\vspace {1mm}
 
 (c) If $X$ is homeomorphic to the $3$-ball, one could also  use  a (finer than Llarull's)  version of the  hypersphericity inequality  due to  the Marques  and Neves, that is  \vspace {1mm}
 
 \hspace {-6mm}{\it a  bound on areas of spherical  minimal surfaces in $X_{\delta \delta}$  with Morse indices  one.}\footnote {This is  theorem 4.9 in  [MN 2011],  which  implies a sharp(!) bound on {\it the  2-waist} of $X\cup_YX$,  which is significantly more precise than  Llarull's  inequality, where, conceivably,   there is  a version of this bound applicable to non-spherical   3-manifolds.}\vspace {1mm}

 (d) A suitable form of the $L_2$-index theorem (see [G 1996])  allows a generalisations of Goette-Semmelmann' theorem to manifolds $Z$ the universal covering of which, rather than $Z$ themselves, are spin.  Accordingly, the inequalities  $[mean]_{Sc\geq 0}$ and  $[mean]_{Sc< 0}$ generalise to manifolds $X$  the universal covering of which are spin.

And, in a similar fashion, these inequalities     generalise to complete non-compact manifolds $X$.


 \vspace{1mm}

  {\it  Warning Exercise.} Find  a counterexample to the following  "proposition" and then find the  mistake in its "proof".
 \vspace {1mm}

 Let $X$ be a closed  orientable Riemannin spin $n$-manifold with 
 $Sc(X)\geq -\varepsilon  $ and let $U\subset X$ be an open connected  subset,  such that the scalar curvature (function) of $X$ is  large on $U$, say
 $$Sc(X, u) \geq 1000n(n-1) \mbox  { for all } u\in U.$$.
 
 "{\sf Proposition}". {\sl  If $\varepsilon>0$,  is sufficiently small, say $\varepsilon \leq 0.001$, then  $X$ admits no area decreasing map 
$f:X \to  S^{n}=S^{n}(1)$  of non-zero degree that sends the complement of $U$ to a point  $s_0\in S^{n}$.}

 \vspace {1mm}
 
 {\it "Proof"}.  Assume otherwise,   let $X'=X\times S^2$ and compose the map  $f'=(f,id): X'\to S^{n}\times S^2$  
 with a smooth {\it distance decreasing} map $S^{n}\times S^2 \to S^{n+2}\left(\frac {1}{10}\right )$ of degree one, which collapses  $\{s_0\}\times S^2$ to a point.   Since the resulting map  $F':X'\to S^{n+2}$ is {\it area decreasing} and since the scalar curvature of $X'$ is {\it everywhere positive}, say for $\varepsilon<1$,  
and since  
$$Sc(X')  > 100(n+1)(n+2)=Sc\left (S^{n+3}\left(\frac {1}{10}\right )\right )$$ 
at {\it all points where the differential of $F'$ doesn't vanish}, namely in $U\times S^2$, the "proof"  follows by the  contradiction with   rough  hypersphericity inequality.
 
\vspace {1mm}
 {\it Hint.} The above argument becomes valid if "no area decreasing map 
$f$" is  replaced by  "no distance  decreasing map $f$".\footnote { The method   of minimal hypersurfaces (which needs no spin)  applied to $U$ allows an alternative proof of this  in most (but not all!)  cases, see [GL 1983] and   [G 2018].}

 \section {Rigidity and Extendability of Metrics with Lower Bounds on the Scalar Curvature.}
 
{\sf To put our  inequalities   $[mean]_{Sc\geq 0}$  and    $[mean]_{Sc< 0}$ to a perspective,  let us remind several earlier results by Miao, Shi-Tam, and Mantoulidis-Miao.
\vspace{1mm}}

 Let $Y_0$  be a  smooth closed   hypersurface in the Euclidean space $ \mathbb R^{n}$,  let $X$ be a compact $n$-manifold with the boundary $Y$ isometric to $Y_0$ and let $mc_0=mc_0(y)$, $y\in Y$, be the mean curvature of $Y_0$ transported to $Y$ by our isometry $Y\leftrightarrow Y_0$.\vspace{1mm}
 
{\huge $\star $}{\large \sf Mean Curvature Rigidity Theorem.} {\it  If $Sc(X)\geq 0$, then  the mean curvature  of $Y \subset X$  can't be greater than that  of $Y_0\subset \mathbb R^{n}$,
$$mean.curv (Y, y)\ngeq mc_0(y),$$ 
unless  $X$ is isometric to the  domain $X_0\subset \mathbb R^{n}$  bounded by $Y_0\subset \mathbb R^{n}$.}\vspace{1mm}
  
  This theorem, which in the case $Y_0=S^{n-1}$ can be regraded  as a (partial)  upgrade of    $[mean]_{Sc\geq 0}$ from extremality to rigidity,  is proved in   [M 2003] and in  [ST 2003]) in three steps, roughly, as follows.\vspace{1mm}

$\star$   Attach $X$ to  the complement $Z_0=\mathbb R^{n}\setminus X_0$ by   the 
 "$\leftrightarrow$"-isometry 
  $$Y=\partial X\leftrightarrow \partial Z =Y_0.$$
$\star$$\star$  Smooth  the   $Y$-corner in the resulting manifold
  $$W=X\cup_Y Z$$
while keeping  $Sc(W)$ everywhere $\geq 0.$\vspace{1mm}

$\star$$\star$$\star$ Use   
{\large  \sf Euclidean Rigidity theorem},  which we have already  met in section 1 and which says in the present notation the following.  \vspace {1mm}

{\Huge $\ast $}  {\sl If  a  complete $C^3$-smooth Riemannin manifold $W$  with  $Sc(W)\geq 0$  is  isometric to $\mathbb R^{n}$
 at infinity,  then $W$ is isometric to $\mathbb R^{n}$.} \vspace {1mm}

 
The following  more satisfactory version of {\huge $\star $} is proven in        [ST 2003]), for convex hypersurfaces. \vspace {1mm}
 
{\huge $\star $} {\huge $\star $} {\large \sf  Integral  Mean Curvature Rigidity Theorem}. {\it  If $Y_0$ is convex, then 
  $$\int_YM(y)dy\leq\int_Y M_0(y)dy,$$
  where the equality $\int M(y)=\int M_0(y)$} implies that $X$ is isometric to $X_0$.\vspace{1mm}
  
{\sc Hyperbolic Remarks.} (\textbf H$_1$) It is shown  in [MM 2016],  Proposition 3.2,  that,  \vspace{1mm}

{\sl Min-Oo's rigidity and positive mass theorems for hyperbolic spaces} 

\hspace {-6mm}  allow  \vspace{2mm}
 
 \hspace {-8mm} {\huge $\star$}\hspace {-1mm}$_{-\kappa }  \hspace{1mm} $   {\it  extensions of   theorems {\huge $\star $} and {\huge   $\star \star $}  to  hypersurfaces $Y_0$ in the hyperbolic spaces $ \mathbb H_{-\kappa}^{n}$ with the sectional curvature $-\kappa$ and manifolds $X$ with $Sc(X)\geq -n(n-1)\kappa$.} \vspace{1mm}
 
 This is can be reformulated in an  especially pleasant  manner  for $n=2$  if  $Y=\partial X$ is homeomorphic to $S^2$ and has $sect.curv(Y) \geq -\kappa$ due to   {\it
 
 Pogorelov's existence (and essential  uniqueness)  theorem for isometric embeddings   of spherical surfaces  with metrics having sectional curvatures $>-\kappa$  to    $\mathbb H_{-\kappa}^{3}$.}\vspace{1mm}

(\textbf H$_2$) The above mentioned results  from    [M 2003], [ST 2003] and [MM 2016]  allow  small perturbations  of the metrics on $Y$ which  can be deformed back   to their original  values  with a minimal loss of positivity of the scalar and the mean curvatures. 

Thus, one easily shows, for instance, that \vspace {1mm}
 
 {\it given a smooth  closed  hypersurface $Y_0$ in $\mathbb H^{n}_{- \kappa}$, there is an {\it effectively computable} $\Delta=\Delta_{N,\kappa} (Y)>0$, such that  if a   metrics $h$ on $Y_0$ is  $C^2$-close to 
 the original  metric  $h_0$ on $Y_0\subset \mathbb H^{n}_{- \kappa}$ induced from $\mathbb H^{n}_{- \kappa}$ 
 $$||h-h_0||_{C^2}\leq \delta<\Delta,$$
 then $Y$ admits no filling by manifolds $X$ with $Sc(X)\geq  -n(n-1)\kappa$ and with the mean curvature of $Y=\partial X$ in $X$ bounded from below by  $M_{\delta}$  for function  $M_\delta$ on $Y_0$  which converges to the mean curvature of $Y_0\subset \mathbb H^{n}_{- \kappa}$ for $\delta\to 0$.}\vspace {1mm}
 
  Whenever applies, this inequality is sharper than  our    $[mean]_{Sc<0}$, but,  unfortunately
    this $\Delta$ is pretty small, something like $1/n^2$,    even for the unit sphere $S^{n-1}\subset \mathbb H^{n}_{- 1}$. \vspace {1mm}
    


(\textbf H$_3$) Min-Oo's proof of his rigidity and the positive mass theorems uses a version of Witten's Dirac operator method, thus, it needs the spin condition.  
Accordingly, the original proofs in [M 2003], [ST 2003] and [MM 2016]  needed $X$ to be {\it spin}, but, as we know now, 
the spin condition is redundant according to [SY 2017]. 

Moreover, 
 Min-Oo's rigidity theorem remains valid for manifolds  $Z= \mathbb H_{-1}^{n}/\Gamma$ for {\it all discrete parabolic} isometry groups,\footnote {An isometry group $\Gamma$ of 
 $H_{-1}^{n}$   is {\it parabolic} if  there is  a horosphere in $H_{-1}^{n}$ invariant under the action of $\Gamma $, or, equivalently,   if    all isometries $\gamma \in \Gamma$ except $id$  keep fixed a {\it unique common  fixed pint} in  the ideal boundary of  $H_{-1}^{n}$. 
 
 Besides,  in the present context,  we also require that the actions of $\Gamma$ on $H_{-1}^{n}$ 
   are {\it free}.  But this, in fact, becomes unnecessary, if instead  of deformations of   $H_{-1}^{n}/\Gamma$, we speak of $\Gamma$-equivariant deformations of   $H_{-1}^{n}$, such that  the $\Gamma$-quotients of the supports of these deformations are compact.} where it reads as follows.  \vspace{1mm}

(*)   {\large \sf    Generalised Min-Oo Rigidity Theorem.} {\sl Parabolic quotient manifolds  $Z=H_{-1}^{n}/\Gamma$ admit no non-trivial "deformations" with compact supports and with $Sc\geq -n(n-1)$},

 where these "deformations" may arbitrarily  change the topology of (a compact region in) $Z$ with {\it no condition on spin}.\vspace{1mm}

 Notice,  that the original Min-Oo rigidity theorem corresponds to the case where $\Gamma=\{id\}$  and that 
 the  rigidity for all $\Gamma$, including Min-Oo rigidity itself,  trivially  follows from the rigidity  of "hyperbolic cusps", i.e. where $\Gamma$ is isomorphic to $\mathbb Z^n$. 
 
Thus,  all of   (*)  reduces to   the case  of cusps,  where the proof  follows by a Schoen-Yau style argument with   a use of minimal $\mu$-bubbles 
   or of minimal surfaces with boundary as it  is  briefly explained in   section 4 below. 

 \section {Manifolds with Negative Scalar Curvatures Bounded From Below.} 
 
{\sf We discuss in  in this section an  approach to the extremality/rigidity problem  based on calculus of variation, that allows sharp inequalities  for $Sc(X)<0$ in some cases.\vspace {1mm}}
 
  Let $X$ be a compact  orientable Riemannin $n$-manifold with, possibly disconnected, boundary, such that
 
$\bullet$  the scalar curvature of $X$ is bounded from below by that of the hyperbolic
 $n$-space $H^{n}_{-1}$,
 $$Sc(X)\geq -n(n-1),$$
 
 $\bullet$ the mean curvature of the  boundary is bounded from below by that of the complement of a horoball in $H^{n}_{-1}$,  
 $$mean.curv(\partial X) \geq -(n-1);$$
 
  $\bullet$ there is a connected  component of  the boundary $\partial X$, call it $Y_+$, the mean curvature of which is  bounded from below  by  $(n-1)$, 
 $$mean.curv(Y_+)\geq n-1.$$

 \vspace{0mm}
 
 {\huge $\bullet$}$_-$ \hspace{1mm}{\sl   If the inclusion homomorphism between   the fundamental groups 
 $$\pi_1(Y_+)\to \pi_1(X)$$
 is injective,
 then   $Y_+$ admits no map with  non-zero degree to  the $(n-1)$-torus   unless  the universal covering of $X$ is isometric to the (infinite)  "band" between too equidistant horospheres in  $\mathbb H^{n}_{-1}$.}
 
   \vspace {1mm}
This is proved in  \S$5\frac {5}{6}$ of [G 1996] for $n\leq 7$ by a "soap bubble"argument which  extends to all $n$  in view of  [SY 2017] and  where the argument from [L 2018] may be also applicable.
 
  \vspace {1mm}

{\it Example and Generalisation.} Let $X_0$ be a compact    {\it hyperbolic}  (i.e. $sect.curv(X_0)=-1$) manifold of dimension $n$ with  {\it totally geodesic} boundary $Y_0$. 

Observe that  the equidistant deformations $Y_d$ $d\geq 0$ of  $Y_0$  in the ambient complete hyperbolic manifold $X_+\supset X$ have $mean.curv(Y_d)\to n-1$ for $d\to \infty$. But, the above implies that  \vspace{1mm}

 {\sf 
no Riemannin manifold $X$ homeomorphic to $X_0$, or just admitting a map with non-zero degree to $X_0$, can have $Sc(X)\geq -n(n-1)$ and 
 $mean.curv(\partial X)
\geq n-1$,} provided,

{\it $Y_0=\partial X_0$, or a finite covering of it admits a non-zero degree map to the $(n-1)$-torus,}

 \vspace{1mm}

This argument doesn't directly apply  if no map $Y_0\to \mathbb T^n $ with $deg\neq 0$ exists.
However, since $X_0$ is {\it isoenlargeable} in the sense of  [G 2018] a 
combination of the arguments from [G 1996]  and  [G 2018]  (with a use of [SY 2017] for $n+1\geq 9$)   shows that

 {\sf the above  conclusion still holds for {\it all  compact     hyperbolic}  manifolds $X_0$  with  {\it totally geodesic} boundaries}.  Namely, 
 $$Sc(X)\geq  -n(n-1) \Rightarrow mean.curv(\partial X)\ngtr n-1.$$
{\sf for all $X$ homeomorphic to $X_0$.}
  \vspace {1mm}
 
 {\it  Application of Symmetrization.}  In many (all?) cases {\huge $\bullet$}$_-$ can be proved by the  torical symmetrization argument from [G 2018]  with no use of "bubbles" at all. A characteristic example is as follows. 
  
Let $X$ be an $n$-dimensional  Riemannin manifold homeomorphic to $Y\times [0,1]$ for $Y=\mathbb T^{n-1}$.\vspace {1mm}   
  
  {\sl If \hspace{0.4mm}$Sc(X)\geq -n(n-1)$ and if the mean curvature of\hspace{0.4mm}
     $ Y_1=Y\times \{1\}\subset \partial X$ is bounded from below by $M_1>n-1$, then 
  $$dist(Y_0=Y\times \{0\}, Y_1) \leq l= \frac {2}{n} \coth^{-1}\frac {M_1}{n-1}, $$}
 {\sf where $\coth^{-1}$ denotes    the inverse function of   
  $\coth x=\frac {e^x+e^{-x}}{e^x-e^{-x}}$.}\vspace {1mm}

In fact (see [G  2018]), the extremal Riemannin  metric $g_{ext}$ with the maximal $l$  is  defined  on the manifold
  $\mathbb T^{n-1}\times (0,l]$,  where it is 
   invariant under the obvious action of the torus $\mathbb T^{n-1}$.
Hence,  $g_{ext}=dt^2+\varphi^2h_{flat}$, where one may assume that $\phi(t)\to 0$ for    $t\to 0$.

Since,
$$mean.curv(g_{ext})= \frac {(n-1)\varphi}{\phi'}$$
by the {\it first variation formula} and  
 $Sc(g_{ext})=-n(n-1)$, 
 {\it Weyl's Riemannian tube formula} shows that   
  $$-n(n-1) = Sc(g_{ext})=-2(n-1)\frac {\varphi''}{\varphi}-(n-1)(n-2)\frac {\varphi'^2}{\varphi^2}.$$

  \vspace {1mm}
 


 
 
    
It follows that the function $f= \frac {\phi}{\phi'}$ satisfies the equation
  $$\frac {f'}{1-f^2} =(\coth^{-1} f)' =\frac {n}{2},$$
   which implies that
   $$f(t)=\coth \frac {tn}{2},$$
 where the domain of definition of this $f$  is  the segment $(0,l]$ for $l= \frac {2}{n} \coth^{-1}\frac {M_1}{n-1}.$
    \vspace {1mm}

  {\sf \large Cuspidal Boundary  Rigidity Theorem.} Let $X$ be  a  complete orientable  Riemannian  $n$-manifold  with compact boundary $Y$, such that  $Sc(X)\geq -n(n-1)$ and $mean.curv(Y)\geq n-1$. \vspace{1mm}
  
  {\huge $\star$}$_-$ \hspace{0.5mm} If  some connected component $Y_0$ of $Y$ admits a map to the $(n-1)$-torus $\mathbb T^{n-1}$ with non-zero degree, which continuously  extends to a  map
  $X\to \mathbb T^{n-1}$,
   then the above argument shows that    \vspace {1mm}
  
\hspace {-7mm} {\sl the universal covering of $X$ is isometric  to a horoball  in the hyperbolic space $\mathbb H^{n}_{-1}$.} \vspace {1mm}
  
(If $X$ is compact and $Y=\partial X$ is connected,  then  maps $Y\to \mathbb T^{n-1}$ with non-zero degrees do not extend to $X$, but this can be sometimes remedied by passing to an infinite coverings $p: \tilde X\to X$, such that the pullbacks 
$p^{-1}(Y)\subset \tilde X$  are  disjoint union of compact manifolds.)
  \vspace {1mm}

{\it  Remark.} The advantage of {\huge $\star$}$_-$ over  the above  {\huge $\bullet$}$_-$ is admission of {\it complete non-compact} manifolds  $X$.  In fact, one may also allow here additional  boundary components with $mean.curv(\geq -(n-1)$, and then {\huge $\star$}$_-$  will imply  {\huge $\bullet$}$_-$.

    \vspace {1mm}

{\it Question.} Is there a   Dirac operator proof of   {\huge $\star$}$_-$ \`a la  Witten and Min-Oo? 

\section {Problems with  $B^2\times \mathbb T^{n-2}$-Manifolds.}  

Let us  discuss  a  possibility of extension of  the above  to a   class of  compact orientable $n$-dimensional  Riemannian manifolds $X$ with boundaries $Y=\partial X$, where the inclusion homomorphism  of the fundamental groups  is non-injective.    

Namely, 
 \vspace{1mm}
 
 (*$_{n}$)  {\sl Let   an $X$ admit a map with non-zero degree to $B^2\times \mathbb T^{n-2}$, where $B^2$ denotes the  disc. \vspace{1mm}
  
  If  $Sc(X)\geq n(n-1) $  and $mean.curv(Y)\geq M >n-1 $, than, {\sf \color {magenta!80 !black} conjecturally,}    the lengths $L $ of the closed curves in $Y$, which  are  non-homologous to zero in $Y$ but 
 homologous to zero in $X\supset Y$, are bounded by the length of the  circle $Y_0$ in the hyperbolic plane $\mathbb H^2_{-1}$ with    $mean.curv(Y_0)=M_0 =\frac {M}{n-1}$.}\vspace{1mm}
 
If $n=2$, this is  a  baby version of the 
 Miao-Shi-Tam rigidity theorem which reads as follows.  \vspace {1mm}

(*$_2$) {\sl   Let  a  compact connected   surface $X$ with boundary $Y$  has  $Sc(X)\geq -2$ 
  and $mean.curv \geq M>1$. 
   Then $X$ is homeomorphic to the disc and 
  the length  of  $Y$  is bounded by that of the circle $Y_0$ in the hyperbolic plane $\mathbb H^2_{-1}$ with    $Mncurv(Y_0)=M$.} \vspace {1mm}
  
  In fact, this is seen by attaching $X$ to the  complement of the ball in   $\mathbb H^2_{-1}$ with the boundary length equal the length of  $Y$   and using the (obvious) rigidity of of the hyperbolic metric  under   curvature increasing  deformations, which have   compact supports in $\mathbb H^2_{-1}$.   
    \vspace {1mm}
 
{\it Exercise.}  State and prove the "dual" inequality for metrics with sectional curvatures bounded from above. \vspace {2mm}

 If   $n\geq 3$, then some (very limited) results can be derived by Miao-Shi-Tam gluing argument  combined with  the rigidity of the manifold 
 $Z=\mathbb H^{n}_{-1}/\mathbb Z^{n-2}$ that  was  stated in (\textbf H$_3$) in section 3  \vspace {1mm}
 
{\it Example.}  Think of  $\mathbb H^{n}_{-1}/\mathbb Z^{n-2}$ as $\mathbb H^2_{-1}\times \mathbb T^{n-2}$
 with a warped (in particular, $ \mathbb T^{n-2}$-invariant) metric and let $X_0 = B^2_0\times \mathbb T^{n-2}$ for some disc $B_0^2\subset \mathbb   H^2_{-1}$. 
 
 Notice that   the boundary $Y_0=\partial X_0$ is the $n$-torus $(\partial  B^2_0=S^1) \times \mathbb T^{n-2}$
   with a warped product metric where  the  certain warping function on the closed  curve   $S^1=\partial B^2_0$ depends on the position and shape of this curve  in  the hyperbolic plane $\mathbb H^2_{-1}$. 
   
   Then, by  the  rigidity of $Z=\mathbb H^{n}/\mathbb Z^{n-2}$,  \vspace {1mm}

   {\sl no Riemannian manifold with the boundary isometric to 
   $Y_0$ can have the scalar curvature $\geq- n(n-1)$ and the mean curvature of the boundary greater than that
   of the original $mean.curv(Y_0\subset Z)$. }  \vspace {1mm}
  
  Also notice that here, similarly to  what was  indicated in   (\textbf H$_3$) in section 3, one may allow  a $C^2$-small perturbation  $h_{new} $  of the original metric $h_0$ in $Y_0$ induced from $Z$. Namely,  \vspace {1mm}

{\sl no Riemannian manifold with the boundary isometric to 
   $(Y_0, h_{new})$ can have the scalar curvature $\geq- n(n-1)$ and the mean curvature of the boundary  greater than that
   of the original $mean.curv(Y_0\subset Z)$ plus one. } \vspace {2mm}

 Our formulation of (*$_2$) as well as  distinguishing   the case of $\mathbb T^{n-1}$-invariant torical domains is motivated by the possibility of  symmetrisation of an arbitrary  $X$ with no decrease  of its  scalar curvature and of  the mean curvature of the boundary  $\partial X$.
 
 Recall (see [G 2018]),  that "the ultimately"   symmetric manifold $Z_0$ carries  the (warped product )  metric $ dt^2+(\sinh^2 t)\cdot h_{flat}$, $0<t<\infty$, which has {\it constant sectional} curvature  {\it only} for $n=2$,  where, if $n=2$, this metric remains non-singular at $t=0$  if   $h_{flat}$  is   the metric of  the circle $S^1_{2\pi}$ of length (exactly)  $2\pi$.
 
 In general, we assign   $h_{flat}$ to  the torus $\mathbb T^{n-1}\times S^1_{2\pi}$ and think of 
 $Z_0$ as a torical warped product over the hyperbolic plane minus a point.
 
 But we don't know 
 
 {\sl whether 
 such metrics are  rigid with respect to deformations with compact (bounded?) supports and $Sc\geq -n(n-1)$.}\vspace {1mm}

 Another, even more annoying, problem  is that the symmetrization in the present case terminates  at a warped product  metric $g_0$ on $X_0=\Sigma^2 \times \mathbb T^{n-1}$, where $\Sigma^2$ is a disc with a Riemannin  metric $\underline g_0$, where we don't control 
  either   $\underline g_0$ or the warping factor that is a function on $\Sigma^2.$
 
 Ideally, we would like the warping factor $w(\sigma)$, $\sigma\in \Sigma^2$,  of this $g_0$ on the (circular)  boundary of  
 $\Sigma^2$  to match the warping function on $\partial B^2_0\subset \mathbb H^2_{-1}$ in   the above example.
 
 But, apparently,  $w(\sigma)$ for  $\sigma \in \partial \Sigma^2$ depends on the geometry of all of $X$, not only on  $\partial X$.
This  prevents us  from a proof of an  even non-sharp version of (*$_2$) by the gluing argument  used for $n=2$.
 
On the other hand one may think of an alternative geometric approach, which, in particular, 
 would deliver   a proof  for the case of  $n=2$   
       by analysing the intrinsic geometry of  the surface  $X$  from  (*$_2$)  rather than by  attaching  something to it.\footnote{A  non-sharp inequality can be derived  from   "strong generalised  concavity"  of the distance function $dist_X (s_1,s_2)$ for $s_1,s_2\in \partial X$, as in the "proof" of unproven corollary in section 9.} 
 \vspace {1mm}
 
{\it  On Surfaces in  3-Manifolds.} If $n=3$, then  additional possibilities are opened by  the existence of isometric imbeddings of   "many" non  rotationally symmetric Riemannin metric on  2-tori to $ \mathbb H^3/\mathbb Z$.
 (Probably, the space of  the  embeddable  metrics  has  codimension one in the space of all metrics.)  
 
 The situation is more satisfactory for imbedding into   $ \mathbb H^3/\mathbb Z^2$, where instead of an ambiguous  "many" one  can definitely say  "all".
 
 In fact,  let $Y =(\mathbb T^2,h)$ be a 2-torus with a  Riemannin metric  where  
 $sectcurv(h) >-1$.
 Then,  by the {\it torical  version of Pogorelov's isometric  embedding theorem}, \vspace {1mm}
 
 {\it there exists   a hyperbolic cusp with the sectional curvature $-1$ and an (essentially unique) isometric embedding $Y\to Z$.}\footnote {The set of   the isometry classes of these "cusps" that are the  quotient manifolds    $ \mathbb H_{-1} ^3/ \Gamma$, $\Gamma=\mathbb Z^2$, is  naturally parametrised by the modular curve of conformal structures on 
 $\mathbb T^2$, where the modular parameter of $Z$, which receives an   isometric embedding of  $\mathbb (T^2, h)$, depends on $h$.} 
 
(This can be exploited similarly to how that was done in [MM 2016]   and  mentioned in 
{\huge $\star$}\hspace {-1mm}$_-\kappa$   in section 4.) 
  \vspace {2mm}
  
  \section {Manifolds  with Corners}
    Let us briefly discuss here what can be done about Sormani's question  {\sf \textbf C } from the introduction.
\vspace {1mm}

Basic examples of Riemannian $n$-manifolds with corners are  (convex) polyhedra in 
$\mathbb R^{n}$ with smooth Riemannian metrics on them.  

In general, a corner structure $P$ on a smooth manifold $X$ with boundary $Y=\partial X$ is given by {\it the   shadow} $\underline P$ pf $P$, that  is a partition of $Y$ into locally closed submanifolds corresponding to  the actual faces of $P$.

For example, the shadow of the corner structure of the cube $[-1, 1]^{n+1}$ on the unit $n$-sphere $S^{n-1}\subset \mathbb R^{n}$ can be seen by radially projecting  the boundary of the cube to this sphere.

More generally, let $X$ be a smooth $n$-manifold with boundary  $Y$  and let $f:Y\to S^{n-1}$ 
be a smooth map  which is transversal to all faces of such a  shadow $\underline P$ in $S^{n-1}$. Then the pullbacks of these faces define a shadow of a corner structure on $X$, where the  corresponding 
corner structure on $X$ is called {\it cubical } if  $X$ is orientable and the map $f$ has {\it non-zero degree.}
 
 ($\square_\circ $)  The simplest instance of this, that we shall use below, is where the map $f$  is a diffeomorphism from the interior of a ball $B\subset Y$ onto $S^{n-1}$ minus the south pole $s_\star\in S^{n-1}$ and where
 all of the complement of $B$ goes to $s_\star$.

\vspace {1mm}



Besides the  scalar  curvature of $X$
and the mean curvatures of its $(n-1)$-faces  an essential geometry of a Riemannin manifold with corners   is carried by  the {\it dihedral angles} $\alpha$ at the   {\it "edges"}  that are the $(n-2)$-faces of $X$, where the difference $\pi-\alpha$ plays  the role of the mean curvature.\vspace{1mm}

An essential  motivation of what we try to do  in this section is the  following. \vspace {1mm}

{\large   \sf Hyperbolic Subrectangular Theorem}. Let $X$ be a compact cubical Riemannian $n$-manifold, where \vspace {1mm}

{\sf (i)  the dihedral angles are  $\leq\frac {\pi}{2}$;

(ii) all  $(n-1)$-faces but  one have positive mean curvatures;

(iii) the exceptional  face have $mean.curv\geq- (n-1)$;

(iv)  the opposite face, which is also called "exceptional",  has $mean.curv>n-1$.}\vspace {1mm}

{\it Then the   scalar curvature of $X$ satisfies:
$$\inf_{x\in X} Sc(X,x)< -n(n-1).$$}

{\it Idea of the proof.} Reflect $X$ around  $2(n-1)$  non-exceptional faces, smooth the resulting $C^0$-metric corners and apply 
{\huge $ \bullet_-$} from section 6 to the resulting "sub-hyperbolic  band".

 In truth, however, such a  smoothing is a mess;  a technically less demanding approach is explained in [G  2014] and in section 11.10 of [G 2018].

(This proof of the theorem is immediate   for  $n=2$, where it is seen  by looking  at minimal geodesic segments between two exceptional faces).  \vspace{1mm}

{\it Rigidity Question.} If {\sf (iv)} is relaxed to  $Mncurv\geq n$, and if $Sc(X,x)\geq -n(n-1)$, then, most likely,   $X$ is isometric to a {\it parabolic rectangular solid} in 
$\mathbb H_{-1}^{n}$, that is defined 
in the {\it horospherical coordinates}  $ x_1,...x_n$ by the inequalities 
$$0\leq x_i  \leq a_i,\hspace {1mm}  i=1,...,n,$$
where, recall, the hyperbolic metric $g$ in these coordinates is  
$$g=dx_1^2 +e^{2x_1}\sum_{i=2}^ndx_i^2.$$ 

We    fail short of directly proving this  because of,   a priori possible, presence of    singularities of the boundaries of  extremal  $X$, which is due to a use  of  the Schoen-Yau style variational argument in an essential step of our proof.
 
 On the other hand, a suitable adapted  {\it Kazdan-Warner perturbation argument}  would, probably,  reduce  the rigidity problem for $Sc(X)\geq -n(n-1) $ to that, where all faces but one are convex and   $Ricci(X) \geq -(n-1)g$; then  the proof  would follow by  Weyl's tube formula.

\vspace {1 mm}

{\it \textbf {Unproven Corollary.}}
Let $X$ be a compact Riemannian $n$-manifold with boundary $Y=\partial X$,   such that 
$$Sc(X)\geq -1\mbox { and } mean.curv(Y)\geq -1.$$
Let $B=B_y(r)\subset Y$, $y\in Y$, be a ball of radius $r$, such that

 the sectional curvatures of $Y$ in this ball are  bounded in absolute values by a constant $\kappa>0$;

the exponential map $T_y(Y)\to Y$ is one-to-one  on    the ball   $B_0(r)\in T_y(Y)$ to $B.$
 
 Then \vspace{1mm}

 {\sl  the infimum of the mean  curvature of $Y$ in the ball $B$ satisfies 
$$\inf_{y\in B} mean.curv(Y,y)\leq {\mathscr M}_n(\kappa +r^{-1})\leqno{[\square_{-1}]}$$
for some (possibly very large) universal continuous function  ${\mathscr M}_n$}.\vspace {1mm}

{\it How we Want to Prove it.}  
The {\it actual proof}   for  $n=2$ is easy.
In fact,  let $X$ be a  surface   with circular boundary $Y=\partial X$, and let $Y_+ \subset Y$ 
be a segment,  such that the following conditions are satisfied.\vspace {1mm}

{\sf (1)  sect.curv$(X)\geq-$1; 

(2) the curvature of the  segment $Y_+$ is $\geq M_+=1+\varepsilon_+ $, $\varepsilon_+>0$;

(3) the curvature of $Y$ in the complement of $Y_+$ is $\geq -1$.}\vspace {1mm}

{\sl Then the length of $Y_+$ is bounded by} 1 000 {\sl  times the length of the circle $Y_{M_+}\subset \mathbb H^2_{-1}$ with the  curvature  $curv(Y_{M_+})=M_+$}.
\vspace {1mm}

{\it Proof.}   Assume, this is easy to justify,   that   the curvature of $Y$ is constant on $Y_+$  and  let  us endow the cylinder $Y\times \mathbb R_+$ to $X$ with a (canonical)  
metric with sectional curvature $ -1$, and such that the boundary  $Y\times \{0\}$ of this cylinder  admits an  isometry, i.e. a  length preserving map 
$$Y\times \{0\}   \leftrightarrow Y,$$
 where this 
isometry  matches  the curvature of these curves, which makes the   resulting metric on the extended manifold 
$$X_+=X\cup_YY\times \mathbb R_+ $$
$C^1$-smooth. (This is essential only on $Y_+$.)

If the length of the segment $Y_+$, now positioned in  $Y\times \mathbb R_+$,  were sufficiently long  then, by an elementary argument, there would exist a 
  deformation of   $Y$ in   $Y\times \mathbb R_+$ to a curve 
$$Y'\subset Y\times \mathbb R_+  \subset X_+,$$ such that

$\bullet$ there are  four 90$^\circ$ corners  on this curve $Y'$.

$\bullet$ the curvatures of one pair of disjoint segments bounded by the corner points  are
$\geq 0$, while the curvatures of the second pair of segments, call them $Y'_+$ and $Y'_-$ are bounded by 
$$curv(Y'_+)\geq 1 +0,01 \varepsilon_+\mbox { and }   curv(Y'_+)\geq -1.$$

This  contradicts    the (obvious in this case)  the above  (sub)rectangular  theorem, applied to the domain $X'$ bounded  by the curve $Y'$  in $X_+$, and thus, the proof is concluded. \vspace {1mm}

{\it Questions.} (i) Which    functions  can be realised as curvatures of  boundaries $Y=\partial X$ of surfaces  $X$, where  $sect.curv(X)\geq \kappa$? 

To get an idea of the expected  patterns of segments of $Y$ with   large and small  curvatures , 
let $Y\subset \mathbb R^2$ be a simple  closed curve, which contains a segment $Y_+$  of length $L_+$, where the curvature  of $Y$ is $\geq 1$. 

It is not hard to show that $Y$ must contain  an open subset   $Y_-$ (possibly, with arbitrarily  many connected components)  of length $L_-$, such that 
$L_-\geq (1-\varepsilon) L_+  $  and 
$$\inf curv(Y_-)=\inf_{y\in Y_-}curv(Y,y) \leq - 1+ \delta, $$
where $\varepsilon, \delta\to 0$ for $L_+\to \infty.$

Probably, all surfaces $X$  with  $sectcurv\geq \kappa$ display similar patterns  of distributions of curvatures of their  boundaries $Y$ for  large $length(Y)$. \vspace{1mm}

(ii) What would be   analogues of the above for (distribution of) the mean curvatures of boundaries $Y$
of $n$-manifolds $X$  with lower bounds on  the sectional,  Ricci  or   scalar curvatures of $X$?

(The case of bounded domains $X\subset \mathbb R^{n}$, $n\geq 3$, already  seems interesting.)

\vspace {2mm}

{\it Idea of a Possible  Proof of $[\square_{-1}]$}.  We want to extend  the above "cornering" argument  to  $n\geq 3$ and create a (sub)rectangular corner structure on $X$, the   shadow of which  will be  induced by a map $B\to S^{n-1}$ as in the above  ($\square_\circ $).

It seems manageable  to bend $Y$   along individual edges, that are the  $(n-2)$-faces of  $X$ (compare with the constructions in  section 11.3, 11.4 in [G 2018]), but it  less clear
 how to do it (I think it is unpleasant but  possible)  at the corners, where the edges meet. 
 
 However, similarly how this was bypassed  in [G 2014[ and in [G 2018]  one can, probably, avoid this problem by 
 performing the  constructions   of corners, reflections and smoothings {\it interchangeably. }
 
For instance, let $X$,  topologically, be  the 3-ball  and let us create two circular corners on its (spherical) boundary  corresponding to the (circular)  boundaries of two non-exceptional faces. 

Then double the resulting $X'$ over this pair, smooth it and arrive at $X_1$, which is now homeomorphic to $B^2\times S^1$, where the shadow of the  corner structure is induced from that on the disc $B^2$.

Proceed as before, now  with the remaining pair of non-exceptional faces and, thus,  arrive at 
$X_2$, which is  homeomorphic to $S^1\times S^1\times [0,1]$  and  to which {\huge $\bullet_{-}$} from section 6 applies. \vspace {1mm}

Conceivably,  the  proof of $[\square_{-1}]$ can be achieved    along these lines. But  artificiality of the argument and the issuing non-sharpness of the result  leave the problem of the correct (re)formulation and of the proof of $[\square_{-1}]$ open.\vspace{1mm}

 On the other hand, the geometry behind this may be interesting in its own right as  suggested by following. \vspace{1mm}

{\large \sf  Question.} Let $X= (X,g)$ be a  Riemannin manifold diffeomorphic to the   $ n$-ball
and let 
 the scalar curvature of $X$ in a neighbourhood $U\subset X$  of $Y=\partial X$ be positive.  Let $V$ be a closed smooth $n$-manifold.\vspace{1mm}

 {\it A packing  with $Sc>0$} of $V$ by (copies of)  $X$ is a Riemannian metric $g_+$ on $V$, along with   isometric embeddings of several  copies of $(X,g)$ to $(V, g_+)$,  such the metric $g_+$ has 
 
\hspace{-5mm} {\it positive scalar curvature in the complement of the images of these embeddings}.

\vspace {1mm}

{\sf Under what conditions on  the intrinsic geometry of $Y=\partial X$, on  the mean curvature of $Y$  and on the topology of $V,$ can $V$ be packed with $Sc>0$  by $X$?}

\vspace{1mm}

More specifically, let $locgeo(Y)$ denote the maximum of two invariants of $Y$
$$locgeo(Y)=\max \left (+\sqrt {|sect.curv(V)|}, \frac{1}{injrad(V)}\right),$$
let    $  top(V)$ means "topology of $V$"  and let ${\mathscr M}_{top(V)}(\ast)$ be a (large) continuous function. \vspace{1mm}

 {\sf For which $V$,  is the inequality $mean.curv(Y)\geq {\mathscr M}_{top(V)}(locgeo(Y))$
sufficient for the existence of such packing?}

\vspace {1mm}

Notice in this regard, that our cornering argument, if it works,  implies  the existence of  such packings  of  the $n$-torus,  provided  ${\mathscr M}={\mathscr M}_{\mathbb T^{n}}$ is  taken sufficiently  large.\footnote{The packings associated with  such "cornering", may have, if you wish,  the complements to the balls  
contained in the $\varepsilon$-neighbourhood of these balls for   arbitrarily small $\varepsilon>0$.}

Also,  such packings   are likely to exist in  other "reflection manifolds".\vspace {1mm}

Contrariwise,  the   largeness of the  curvature of  the boundary curves  serves as an obstruction for  packing  with $Sc>0$ of 2-spheres  by discs.
 
More  generally,  {\it aspherical}\footnote {A topological space is called aspherical if its universal covering is contractible.} manifolds $V$ of all dimensions  have a chance  to be packable, (more realistically, {\it approximately} packable in some  sense),   under similar conditions    but it is rather improbable  for the spheres, even if "$locgeo$" is replaced by a more comprehensive global invariant of $X$. 
 But finding relevant  examples of   non-extendability of metrics with $Sc>0$   seems difficult.

(All what matters  for the existence of  an  extension with $Sc>0$  is the geometry of $X$ near  $Y=\partial X$. But  what happens in $X$ far from $Y$ may harbour  obstructions to such extensions.)

\vspace{2mm}

{\large \sf Bound on Hight of Fat Cylinders.} Let us indicate another application of the cornering construction,
where there is an instance where   no higher order corners difficulty  appears.
 \vspace{1mm}

Let $X$  be a compact connected orientable  $n$-manifold with boundary $Y$ and let $Y_{bot}\subset $ and $Y_{top}$ be two disjoint smooth connected domains in $Y=\partial X$ called the {\it bottom} and the {\it  top.} \vspace{1mm}
 
 The basic example is where $X$ is  the cylinder, $X_0= B^{n-1}\times[0,1]$;   in general, 
 $X=(X, Y_{bot} ,X_{top})$ is called   "cylinder"  if it admits a continuous map  
  $f: X\to B^{n-1}\times[0,1]$,   such that $\partial X$ is sent to $\partial  X_0$,  
 $$f^{-1} (B^{n-1}\times \{0\}) =Y_{bot},\mbox { \hspace {1mm} }f^{-1} (B^{n-1}\times \{1\})=Y_{top}$$
 and $deg(f)\neq 0$. \vspace{1mm}
 
 Denote by $Y_\varocircle\subset Y$   the "side" of  the "cylinder" $X$, that is the pullback
 $f^{-1}(S^{n-2}\times[0,1])$ for  the boundary sphere $S^{n-2}=\partial B^{n-1}$ and let the following conditions be satisfied.

$\bullet $    $Sc(X)\geq n(n-1)$,
  
 $\bullet $  $locgeo(Y_\varocircle) \leq1$,
  
  $\bullet $  $mean.curv(Y_\varocircle)\geq \underline M$.\vspace {1mm}
  
 {\it \textbf {Second Unproven Corollary.}}  {\it The "hight" of $X$ is bounded by a universal function of $\underline M$:
 $$ dist_X(Y_{bot}, Y_{top})\leq {\mathscr H}_n(\underline M),$$ 
  such that 
  $${\mathscr H}_n(\underline M)\to \frac {2\pi}{n}\mbox {  for $\underline M\to \infty$}.$$}
 
{\it Idea  of a Possible  Proof.}  If $n=3$, bend $Y_\circ$ along four segments between 
the top and the bottom of $X$ and, thus,  produce  four edges with dihedral angles $\leq \frac {\pi}{2}$.
This seems non-difficult.  

 More problematic is to create similar $n$-cubical structure for $n\geq 4$
along the lines indicated  in the above  "possible proof" of  $[\square_{-1}]$.

 But if this is accepted, then the {\sl sub-rectangular $\frac{2\pi}{n}$-inequality}  from section 11.10 in [G  2018]  could be applied and the proof  would follow.

\vspace{1mm}

{\it Remarks/Questions.}  (a) The apparent examples suggest that 
 $$ dist_X(Y_{bot}, Y_{top})\to 0 \mbox {  for $\underline M\to \infty$},$$
 where it is, indeed, so for $n=2$  by an elementary argument.
 
 But it is unclear whether this is true or false for $n\geq 3$.

\vspace{1mm}

(b) Can one replace  the bound on $locgeo(Y_\varocircle)$ by a lower bound on some kind of  size of $Y_\varocircle$?

Such a "size"  can be conveniently defined with the above map $f$ in the definition  of the "  cylinder"
structure on $X$  as follows.

Restrict this $f$ to $Y_\varocircle$, observe that it sends 
$Y_\varocircle\to S^{n-2}\times [0,1]$, compose this  with the projection $ S^{n-2}\times [0,1] \to S^{n-2}$ and denote
the resulting map  by
 $\underline f: Y_\varocircle\to S^{n-2}$.
 
  If $n\geq 3$, denote by $\underline Rad_{S_{n-2}}(Y_\varocircle)$ the supremum of the numbers $R$ such that $
\underline f$ is homotopic to a smooth  map $\underline f': Y_\varocircle\to S^{n-2}$, where the norm of the differential of $\underline f'$ is bounded by
$$ \sup_{y\in Y_\varocircle} ||d\underline f'(y)||\leq \frac {1}{R}.$$
And If $n =2$, let  $\underline Rad_{S_{n-1}}(Y_\varocircle)=1$ for all $X$.\vspace {1mm}

 {\sc Conjecture.} Let $Sc(X)\geq 0$  and $ mean.curv(Y)\geq \underline M$.  If the product  
 $$\Theta(X)=\underline M\cdot \underline Rad_{S_{n-2}}(Y_\varocircle)$$ (this $\Theta$ speaks for  "fatness" of $X$) is sufficiently large, say 
 $$\Theta(X)> \Theta_{crit},$$
   then 
$$ dist_X(Y_{bot}, Y_{top})\leq \lambda_n (\Theta(X)) $$
for some continuous monotone decreasing  function $\lambda_n(\Theta)$
defined for $\Theta>  \Theta_{crit}$. \vspace {2mm}

The best one may expect here is that  $$\Theta_{crit}=n-2$$ 
and that 
$$\lambda_n (\Theta)\to 0 \mbox { for }\Theta\to \infty.$$
But any bound $\Theta_{crit}\leq const_n$  will be welcome and the inequality
$$\lim_{\Theta\to \infty}\lambda_n (\Theta) \leq \frac {2\pi}{n-1}$$
will be quite satisfactory. \vspace {1mm}

And  it is conceivable that    all of the above  holds for $Sc(X)\geq \sigma,$ for all  $-\infty <\sigma<\infty$,
with  $\Theta_{crit}= \Theta_{crit}(\sigma)$ and  $\lambda_n=\lambda_n(\Theta , \sigma).$

 \section {Further Conjectures and Problems.}

 Probably,  all we know (and don't know) about manifolds with boundaries extends to manifolds with corners, albeit  by   no means automatically.  
In fact, it will be more productive to study manifolds with corners for their own sake rather than as intermediates in the arguments  concerning smooth manifolds.

For instance,   problem {\large {\sf \textbf  A}}  from the introduction becomes even more interesting for manifolds with corners where it reads as follows.

\vspace {1mm}

 Let $X$ be a  compact  $n$-dimensional Riemannin manifold with corners and $Y_0$ be an $(n-1)$-dimensional face of $X$. \vspace {1mm}

{\sl Find an upper bound on the mean curvature of this face, in terms of  
lower bounds on the scalar curvature of $X$ and the mean curvatures of the remaining $(n-1)$-faces of  $X$  and  the dihedral angles along   $(n-2)$-faces,  as well the some geometric invariants of the face $Y_0$, e.g. its "size" and its own   scalar  curvature.  }\vspace {1mm}

The topologically simplest  instances of this are as follows.\vspace {1mm}

$\bullet $ The half-ball $B_+^{n}$, that is the intersection of the ball with the halfspace $x_1\geq 0$ where the $n$-ball $B_+^{n}\cap \mathbb R^{n-1}$ is taken for $Y_0$;

$\bullet $  the {\it  cylinder}
 $X= B^{n-1} \times [0,1]$,    where $B^{n-1}\times \{0\}$ is taken for $Y_0$;

$\bullet $ products of smooth manifolds with boundaries, e.g. cubes $[0,1]^{n}$
and close relatives of cubes --  $n$-diamonds and  $n$-simplices.\vspace {1mm}

But it is more challenging to  understand the geometry of $X$ when the combinatorics of the corner structure becomes  more complicated. \vspace {1mm}

For instance, define $CombRad_{n-1}(Y)$, for  a closed orientable $(n-1)$-manifold   $Y$ partitioned  into faces,  as the supremum of the numbers $R$, such that 
$Y$ admits a continuous  map of non-zero degree to the sphere $S^{n-1}(R)$, such that 
the images of all faces  have diameters at most one.
\vspace {1mm}

{\sc {\Large $\hexagon$} Conjecture.}    {\sf Let $X$ be a compact orientable Riemannin $n$-manifold with corners, such that  $Sc(X)\geq 0$, and all ${n-1}$-faces have positive  mean curvatures. Let
  $\overline {\mathscr A}(X)$ denote the  supremum of the dihedral angles of $X$ taken over all points of all edges ($(n-2)$-faces)  of $X$. }
  
  {\it Then 
$$\pi -\overline {\mathscr A}(X)\leq \frac {const_n}{CombRad_{n-1}(\partial X)} $$
for some universal (possibly large) constant $const_n$.}
\vspace {1mm}

{\it Admission.}  I haven't check this even in the case  of convex polyhedra 
$X\subset \mathbb R^{n}.$

\vspace {2mm}

\vspace {1mm}
   
 Below is  another kind of   conjecture,  where the difficulty resides in the gap between a  use of  Dirac theoretic methods and those of  minimal hypersurfaces.   \vspace {1mm}
 
 {\sc  {\large $ \bigcirc\hspace {-1.8mm}=\hspace {-0.8mm}=$} Conjecture.}  {\sf There exists a (possibly very large) constant $\sigma_n>0$ with the following property. }\vspace {1mm}

{\it If    a compact orientable  Riemannian $n$-manifold with boundary $Y$ has  $Sc(X)\geq \sigma_n$, then every smooth  {\it area decreasing} map  $f$ from $X$ to the  unit $n$-sphere, which {\it sends 
the $1$-neighbourhood  (unit collar)  of $Y\subset X$ to a single point},  has degree $deg(f)=0$.}

    \vspace {1mm}

We conclude by  formulating an instance  of a general {\it stability problem  for $Sc\geq \sigma$},    
in the spirit of [S 2016], which may have  a satisfactory  solution, at least in dimension 3. \vspace {1mm}

 {\sc Spherical Stability  Problem.} {\sf Describe the geometry of   closed Riemannian orientable  $n$-dimensional  manifolds $X_i$, such that  \vspace {1mm}
 
 \hspace {20mm}$\inf Sc(X_i) \to n(n-1)$ and $Rad _{S^{n}}(X_i)\to 1$.} \vspace {1mm}
    
One expects here, in view of Llarull's theorem, that  manifold $X_i$  with  $Sc(X_i) \geq  n(n-1)-\varepsilon_i$
and $Rad _{S^n}(X_i)\geq 1 -\varepsilon_i $  must look, approximately,  as the unit sphere  $S^n$ with an extra staff attached by "$\varepsilon_i$-narrow bridges" to it. 

 In other words,  certain equidimensional submanifolds $U_i\subset X_i$  with (small?) boundaries must  somehow  converge (sub-converge?),  e.g. in  {\it the intrinsic flat topology} (see [S 2016]), to $S^n$.
 \vspace {1mm}

 Besides $S^n$, the stability problem arises for all,   not only spherical, {\it length extremal} and 
{\it area extremal} closed manifolds with positive scalar curvatures  as well as    for extremal manifolds
 with boundaries (as in (d) of section 1), where the situation is even less clear.\vspace {1mm}

It   seems  helpful for   developing  an idea of what happens  to scalar curvature in this regard, up to a point of     formulating   conjectures,  to look at the corresponding problem(s) for hypersurfaces with positive mean curvatures. 
 
 For instance,

 {\sf  what is the geometry of   (limits of  sequences  of) smooth bounded  domains $U_i\subset  \mathbb R^{n+1}$, which contain the unit ball $B^{n+1}_0(1)\subset \mathbb R^{n+1}$  and such that the mean curvatures of the boundaries of these domains are bounded from  below by  $M_i\to n$
  for $i\to \infty$. }

\vspace {1mm}
 
 We conclude by reiterating the problem suggested  by Pengzi Miao:

 {\sf Is there a {\it lower bound on the volumes}  of  manifolds $X $  with $Sc(X)\geq \sigma$, say for $\sigma= -1,0,+1$, in terms of
 
 {\sf intrinsic geometries of their boundaries $Y=\partial X$  and mean curvatures of }$Y$? }
 
    For instance, 
    
    Let a compact Riemannian  $n$-manifold $X$ satisfies: 
    
  {\sf  $\bullet$ \hspace {1mm}   $Rad_{S^{n-1}}(\partial X)\geq 1$;
   
    $\bullet$ \hspace {1mm} $mean.curv (\partial(X))\geq n-1$;
 
  $\bullet$ \hspace {1mm}  $Sc(X) \geq  -\varepsilon $.

    Does  the volume of $X$  relate to that of the unit ball  as follows.   
      $$vol(X)\geq vol (B^n(1))-\delta
      \mbox {  for }\delta=\delta(\varepsilon)\underset {\varepsilon \to 0}\to 0 ?$$}\vspace {1mm}

  This seems  unclear (do I miss something   obvious?) even if $X$ has sectional curvature $\geq -\varepsilon.$
  

 \vspace{7mm}
 


  \vspace{1mm}
 
\hspace {40mm} {\large \sc  \textbf {References}}  \vspace{1mm}  \vspace{3mm}

[BM 2011] S. Brendle and F.C. Marques, {\sl Scalar curvature rigidity of geodesic balls in
$S^n,$} 

\url {https://projecteuclid.org/download/pdf_1/euclid.jdg/1321366355}
\vspace {2mm}

[BMN 2010] S. Brendle, F.C. Marques, A. Neves. {\sl Deformations of the hemisphere that increase scalar curvature,}

\url {https://arxiv.org/abs/1004.3088}
\vspace {2mm} 

[G 1996] M. Gromov,  {\sl  Positive curvature, macroscopic dimension, spectral gaps and higher signatures,}

\url {http://cds.cern.ch/record/284470/files/SCAN-9507071.pdf}

\vspace {2mm}

[G 2014] M. Gromov,  {\sl Dirac and Plateau billiards in domains with corners,}

\url{ https://www.ihes.fr/~gromov/wp-content/uploads/2018/08/Plateauhedra_modified_apr23.pdf}


\vspace {2mm}

[G 2018] M. Gromov, {\sl Metric Inequalities with Scalar Curvature,}

\hspace {-1mm}\url {https://arxiv.org/abs/1710.04655}

\vspace {2mm}

 [GL 1980] M.Gromov, B. Lawson,  {\sl Spin and scalar curvature in the presence of a fundamental group}  Ann. of Math. 111 (1980), 209-230.\vspace {2mm}

[GL 1983] M. Gromov, B. Lawson, {\sl Positive scalar curvature and the Dirac operator on complete Riemannian manifolds} 

\url {http://www.numdam.org/article/PMIHES_1983__58__83_0.pdf}\vspace {2mm}

[GS 2000]  S. Goette, U. Semmelmann,   {\sl Scalar Curvature on Compact Symmetric Spaces},
   J. Differential Geom. Appl.
, 16(1):65-78, 2002.
   
   \url {https://arxiv.org/abs/math/00   10199}

   \vspace {2mm}

[Ll 1998] M. Llarull,     {\sl  Sharp estimates and the Dirac operator,} Mathematische Annalen
1998, Volume 310, Issue 1, pp 55-71. 

\url { https://link.springer.com/article/10.1007/s002080050136 }\vspace {2mm} 

[L 1999] J. Lohkamp  
   Scalar curvature and hammocks,
 
  https://link.springer.com/article/10.1007\%2Fs002080050266\vspace {2mm}

[L 2018]  J. Lohkamp, {\sl  Minimal Smoothings of Area Minimizing Cones},

\url {https://arxiv.org/abs/1810.03157 }

\vspace {2mm}

[M 2003] P. Miao, {\sl Positive Mass Theorem on Manifolds admitting Corners along a Hypersurface,}

\url {https://arxiv.org/abs/math-ph/0212025}

\vspace {2mm}

[MM 2016] 
C. Mantoulidis, P. Miao, {\sl Total mean curvature, scalar curvature, and a variational analog of Brown-York mass,}
\url{ https://arxiv.org/abs/1604.00927}

\vspace {2mm}

[MN 2011]  F. Marques, A. Neves {\sl Rigidity of min-max minimal spheres in three-manifolds,}
\url {https://arxiv.org/pdf/1105.4632.pdf} \vspace {1mm} 

\vspace {2mm} 

[ST 2003]   Y.G. Shi, L.-F. Tam,
{\sl Positive mass theorem and the boundary behaviours of
compact  manifolds  with  nonnegative scalar  curvature,}

https://arxiv.org/abs/math/0301047
\vspace {2mm} 

[S 2016]  C. Sormani, {\sl  Scalar Curvature and Intrinsic Flat Convergence,}

 \url {https://arxiv.org/abs/1606.08949}
\vspace {2mm}

[SY 2017] R. Schoen, S-T. Yau, {\sl Positive Scalar Curvature and Minimal Hypersurface Singularities},
\url{  https://arxiv.org/abs/1704.05490}\vspace {2mm}

[W 1970] A. Weinstein, {\sl Positively curved $n$--manifolds in $R^{n+2}$,}
J. Differential Geom.
    Volume 4, Number 1 (1970), 1-4.

\url {https://projecteuclid.org/euclid.jdg/1214429270}

\end {document}